\documentclass[a4paper,12pt]{article}
\usepackage{}
    \usepackage[top=2.5cm,bottom=2.5cm,left=2.5cm,right=2.5cm]{geometry}
    \usepackage{amsmath, amssymb}
    \allowdisplaybreaks

\usepackage{mathrsfs}
\usepackage{amsmath, epsfig, cite}
\usepackage{amssymb}
\usepackage{amsfonts}
\usepackage{latexsym}
\usepackage{graphicx}

\newtheorem{thm}{Theorem}[section]
\newtheorem{cor}[thm]{Corollary}

\newtheorem{definition}{Definition}[section]

\newcommand{\qed}{{\hfill\rule{4pt}{7pt}}}
\def\pf{\noindent {\it Proof.} }

\numberwithin{equation}{section}

\makeatletter \@addtoreset{equation}{section} \makeatother

\begin{document}
\thispagestyle{empty}

\vspace*{20mm} \pagestyle{empty}
\begin{center}
{\Large\bf Novel Inequalities for Generalized Graph Entropies\\[2mm] Revisited, Graph Energies and Topological Indices\footnote{Supported
by NSFC, ``973" program and PCSIRT.}}
\end{center}

\begin{center}
\small Xueliang Li, Zhongmei Qin, Meiqin Wei\\[2mm]
\small Center for Combinatorics and LPMC-TJKLC, Nankai University, \\ Tianjin 300071, P. R. China.\\[2mm]
\small Ivan Gutman\\
\small Faculty of Science, University of Kragujevac,
P. O. Box 60, 34000 Kragujevac, Serbia\\[2mm]
\small Matthias Dehmer\\
\small Institute for Bioinformatics and Translational Research, UMIT, Hall in Tyrol, Austria  \\
\small Department of Computer Science, Universit\"{a}t der Bundeswehr M\"{u}nchen, Werner-Heisenberg-Weg 39, 85577 Neubiberg, Germany\\[1mm]
\small E-mail: lxl@nankai.edu.cn; qinzhongmei90@163.com; weimeiqin8912@163.com; gutman@kg.ac.rs; matthias.dehmer@umit.at\\[2mm]

\end{center}

\begin{center}
{\bf Abstract}
\end{center}

The entropy of a graph is an information-theoretic quantity which expresses the complexity of a graph \cite{DM1,M}. After Shannon introduced the definition of entropy to information and communication,
many generalizations of the entropy measure have been proposed, such as R\'{e}nyi entropy and Dar\`{o}czy's entropy. In this article,
we prove accurate connections (inequalities) between generalized graph entropies, distinct graph energies and topological indices. Additionally, we obtain some extremal properties of nine generalized
graph entropies by employing distinct graph energies and topological indices.

{\flushleft\bf Keywords}: generalized graph entropies; graph energies; graph indices

\section{Introduction}
~~The entropy of a probability distribution can be interpreted not only as a measure of uncertainty but also as a measure of information. As a matter of fact, the amount of information, which we get when we observe the result of an experiment, can be taken numerically equal to the amount of uncertainty concerning the outcome of the experiment before carrying it out. Shannon first introduced the definition of entropy to information and communication. Moreover, studies of the information content of graphs and networks were initiated in the late $1950$s following publication of the widely cited paper \cite{SW} of Shannon. Later, entropy measures were used to investigate the complexity of graphs associated
with machine learning procedures. And fortunately, such measures have been proved useful to investigate several important properties of a graph. The broad range of research on entropy measures and graphs is exemplified in \cite{B1,B2,M,TCM}. Early contributions in this field inspired researchers in various disciplines to apply entropy measures to the analysis of structures. Various information theoretic measures (and also noninformation-theoretic measures) and other techniques have been developed to determine the structural complexity of molecular structures and complex networks. And an up-to-date review on graph entropy measures has recently been published by Dehmer and Mowshowitz \cite{DM1}.

~~It is worth mentioning that various graph entropy measures have been developed, see \cite{B1,DM1,M}. For example, partitions based on several graph invariants, such as vertices, edges and distances have been used to assign a probability distribution to a graph. In \cite{B3}, Bonchev proposed the magnitude-based information indices, while the topological information content was developed by Rashevsky \cite{RL}. Moreover, so-called generalized graph entropies have been investigated due to Dehmer and Mowshowitz by applying generalized entropy measures, see \cite{DM2}. The innovation represented by these generalized entropy measures is their dependence on the assignment of a probability distribution to a set of elements of a graph. Rather than determine a probability distribution from properties of a graph, one is imposed on the graph independently of its internal structure. Applying the graph energy and the spectral moments, M. Dehmer, X. Li and Y. Shi \cite{DLS} gave accurate connections between graph energy and generalized graph entropies, which were introduced by Dehmer and Mowshowitz \cite{DM2}. Also, some extremal properties of the generalized graph entropies are described and proved in their article \cite{DLS}.

~~In this article, we focus on the mentioned generalized graph entropy measures and express those quantity using distinct graph energies and topological indices.
This article is organized as follows. Section $2$ reviews existing entropies defined on graphs and generalized graph entropies. In Section $3$, we obtain some extremal properties of nine generalized
graph entropies by employing distinct graph energies and topological indices. Moreover, we state some inequalities for
generalized graph entropies. The paper ends with a short summary and conclusion.

\section{Preliminaries}
~~Graph entropies can be divided into two classes. The first class is based on an equivalence relation defined on the set X of elements of a graph while the second, introduced by Dehmer \cite{D},
is not based on partitions induced by equivalence relations. To define these measures, a probability value to each vertex $v_i\in V$ is assigned, and we obtain the following probability distribution
$$
(p^f(v_1),p^f(v_2),\ldots,p^f(v_n)), \quad\quad |V|:=n,
$$
where \cite{D}
$$
p^f(v_i):=\frac{f(v_i)}{\sum\limits_{j=1}^nf(v_j)}
$$
and $f$ is an information function mapping graph elements (e.g., vertices) to the non-negative reals. The entropy of the underlying graph topology is here defined by \cite{D}:
$$
I_f(G):=-\sum\limits_{i=1}^n\frac{f(v_i)}{\sum\limits_{j=1}^nf(v_j)}
\log\frac{f(v_i)}{\sum\limits_{j=1}^nf(v_j)}.
$$

~~Actually, many generalized entropies have been proposed after the seminal paper of Shannon \cite{SW}. Here, we introduce two important examples of entropy measures: R\'{e}nyi entropy \cite{R} and Dar\`{o}czy's entropy \cite{DJ}. The R\'{e}nyi entropy is defined by
$$
I_\alpha^r(P):=\frac{1}{1-\alpha}
\log\left(\sum\limits_{i=1}^n(p(v_i))^\alpha\right), \quad\quad \alpha>0 \  \text{and} \  \alpha\neq 1,
$$
where $P:=(p(v_1),p(v_2),\ldots,p(v_n))$. And the Dar\`{o}czy's entropy is
$$
H_n^\alpha(P):=\frac{\sum\limits_{i=1}^n\left((p_i)^\alpha\right)-1}
{2^{1-\alpha}-1}, \quad\quad \alpha>0 \  \text{and} \  \alpha\neq 1,
$$
where $P:=(p_1,p_2,\ldots,p_n)$. In \cite{DM2}, Dehmer and Mowshowitz introduced a new class of measures (so-called generalized measures here) that derive from functions such as those defined by R\'{e}nyi
entropy, Dar\`{o}czy��s entropy and the quadratic entropy function discussed by Arndt \cite{A}.
\begin{definition}
Let $G$ be a graph of order $n$. Then
$$(i)  \quad \quad ~~~~~~I^1(G)=\sum_{i=1}^n\frac{f(v_i)}{\sum\limits_{j=1}^nf(v_i)}
\left[1-\frac{f(v_i)}{\sum\limits_{j=1}^nf(v_i)}\right];  \quad  \quad  \quad\quad \quad \quad\quad \quad \quad\quad \quad$$

$$(ii)  \quad \quad ~~~~~I^2_\alpha(G)=\frac{1}{1-\alpha}\log\left(\sum_{i=1}^n
\left(\frac{f(v_i)}{\sum\limits_{j=1}^nf(v_i)}\right)^\alpha\right),\quad \alpha \neq 1 ; \quad  \quad\quad \ \quad\quad\quad $$

$$(iii)  \quad \quad \quad\ I^3_\alpha(G)=\frac{\sum\limits_{i=1}^n\Bigg(\frac{f(v_i)}
{\sum\limits_{j=1}^nf(v_i)}\Bigg)^\alpha-1}{2^{1-\alpha}-1},\quad \alpha \neq 1.  \quad  \quad\quad \quad\quad \quad\quad \quad \quad\quad\quad$$

\end{definition}

~~Let $G$ be a graph of order $n$ and $M$ be a matrix related to the graph $G$. Denote $\mu_1, \mu_2, \ldots, \mu_n$ be the eigenvalues of $M$ (or the singular values for incidence matrix). If $f:=|\lambda_i|$, then \cite{DSV}
$$
p^f(v_i)=\frac{|\mu_i|}{\sum\limits_{j=1}^n|\mu_i|}.
$$
Therefore, the generalized graph entropies are defined as follows:
\begin{eqnarray}
&(i)&I^1(G)=\sum_{i=1}^n\frac{|\mu_i|}{\sum\limits_{j=1}^n|\mu_j|}
\left[1-\frac{|\mu_i|}{\sum\limits_{j=1}^n|\mu_j|}\right]; \label{eq1} \\
&(ii)&I^2_\alpha(G)=\frac{1}{1-\alpha}
\log\left(\sum_{i=1}^n\left(\frac{|\mu_i|}
{\sum\limits_{j=1}^n|\mu_j|}\right)^\alpha\right),\quad \alpha \neq 1 ; \label{eq2}\\
&(iii) &I^3_\alpha(G)=\frac{\sum\limits_{i=1}^n\left(\frac{|\mu_i|}
{\sum\limits_{j=1}^n|\mu_j|}\right)^\alpha-1}{2^{1-\alpha}-1},\quad \alpha \neq 1.  \label{eq3}
\end{eqnarray}
Especially, for the first generalized graph entropy $I^1(G)$, we have
\begin{align*}
I^1(G)&=\sum_{i=1}^{n}\frac{|\mu_i|}{\sum\limits_{j=1}^n|\mu_j|}
\left[1-\frac{|\mu_i|}{\sum\limits_{j=1}^n|\mu_j|}\right]\\
&=1-\frac{1}{\left(\sum\limits_{j=1}^n|\mu_j|\right)^2}
\sum\limits_{i=1}^{n}|\mu_i|^2.
\end{align*}

\section{Extremal properties of generalized graph entropies}
~~In this section, we will introduce nine generalized graph entropies of distinct graph matrices. As follows, accurate connections between the entropies and energies or topological indices are proved. Moreover, we examine the extremal properties of the above stated entropies.

1. Let $Q(G)$ be the signless Laplacian matrix of graph $G$. Then $Q(G)=D(G)+A(G)$, where $D(G)={\rm diag}(d_1, d_2, \ldots, d_n)$ denotes the the diagonal matrix of vertex degrees of $G$ and $A(G)$ is the adjacency matrix of $G$. Let $q_1, q_2, \ldots, q_n$ be the eigenvalues of $Q(G)$. Obviously, $q_i \geq 0$, $\sum\limits_{i=1}^n q_i =2m$ and $\sum\limits_{i=1}^n q_i^2=tr(Q^2(G))=\sum\limits_{i=1}^nd_i^2+\sum\limits_{i=1}^nd_i$. Then we have the following theorem.
\begin{thm}\label{thm1}
Let $G$ be a graph with $n$ vertices and $m$ edges. Then for $\alpha \neq 1$, we have
\begin{eqnarray}
&(i)  &I_Q^1(G)=1-\frac{1}{4m^2}(M_1+2m),
\label{eq11}\\
&(ii)&I_\alpha^2(G)=\frac{1}{1-\alpha}\log\frac{M_\alpha^*}{(2m)^\alpha}, \label{eq12}\\
&(iii) &I_\alpha^3(G)=\frac{1}{2^{1-\alpha}-1}
\Bigg(\frac{M_\alpha^*}{(2m)^\alpha}-1\Bigg),\label{eq13}
\end{eqnarray}
where $M_1$ denotes the first Zagreb index and $M_\alpha^*=\sum\limits_{i=1}^n|q_i|^\alpha$.
\end{thm}

\pf By substituting $\sum\limits_{i=1}^n q_i =2m$ and $M_1=\sum\limits_{i=1}^n q_i^2$ into equality (\ref{eq1}), we have
\begin{align*}
I_Q^1(G)&=1-\frac{1}{\left(\sum\limits_{j=1}^n|q_j|\right)^2}
\sum\limits_{i=1}^{n}|q_i|^2\\
&=1-\frac{1}{(2m)^2}\Bigg(\sum_{i=1}^nd_i^2+\sum_{i=1}^nd_i\Bigg)\\
&=1-\frac{1}{4m^2}(M_1+2m).
\end{align*}

The other two equalities can be obtained by substituting $\sum\limits_{i=1}^n q_i =2m$ and $M_\alpha^*=\sum\limits_{i=1}^n|q_i|^\alpha$ into equalities (\ref{eq2}) and (\ref{eq3}), respectively.\qed

Equality (\ref{eq11}) gives the accurate relation between $I_Q^1(G)$ and the first Zagreb index $M_1$. Obviously, for a graph $G$, each upper (lower) bound of the first Zagreb index $M_1$ can deduce a lower (an upper) bound of $I_Q^1(G)$. Moreover, we obtain the following extremal properties of the general graph entropy $I_Q^1(G)$ employing some known bounds in \cite{IIL}.
\begin{cor}
i. For a graph $G$ with $n$ vertices and $m$ edges, we have
$$I_Q^1(G) \leq 1-\frac{1}{2m}-\frac{1}{n}.$$
ii. Let $G$ be a graph with $n$ vertices and $m$ edges. The minimum degree of $G$ is $\delta$ and
the maximum degree of $G$ is $\Delta$. Then
$$I_Q^1(G)\geq 1-\frac{1}{2m}-\frac{1}{2n}-\frac{\Delta^2+\delta^2}{4n\Delta\delta},$$
with equality if and only if $G$ is regular graph, or $G$ is a graph whose vertices have exactly two degrees $\Delta$ and $\delta$ with $\Delta+\delta$ divides $\delta n$ and there are exactly $p = \frac{\delta n}{\delta+\Delta}$ vertices of degree $\Delta$ and $q=\frac{\Delta n}{\delta+\Delta}$ vertices of degree $\delta$.
\end{cor}\qed

2. Let $\mathscr{L}(G)$ and $\mathcal {Q}(G)$ be the normalized Laplacian matrix and the normalized signless Laplacian matrix respectively. By definition, $\mathscr{L}(G)=D(G)^{-\frac{1}{2}}L(G)D(G)^{-\frac{1}{2}}$ and $\mathcal {Q}(G)=D(G)^{-\frac{1}{2}}Q(G)D(G)^{-\frac{1}{2}}$, where $D(G)$ is the diagonal matrix of vertex degrees, and $L(G)=D(G)-A(G)$, $Q(G)=D(G)+A(G)$ are, respectively, the Laplacian and the signless Laplacian matrix of graph $G$. Denote the eigenvalues of $\mathscr{L}(G)$ and $\mathcal {Q}(G)$ by $\mu_1, \mu_2, \ldots, \mu_n$ and $q_1, q_2, \ldots, q_n$ respectively. Then we have $\mu_i \geq 0$, $q_i \geq 0$, $\sum\limits_{i=1}^n \mu_i =\sum\limits_{i=1}^n q_i=n$ and $\sum\limits_{i=1}^n \mu_i^2
=\sum\limits_{i=1}^n q_i^2=n+2 \sum\limits_{i\sim j}\frac{1}{d_id_j}$. As follows, the equality relationship between the generalized graph entropy $I_{\mathscr{L}(\mathcal {Q})}^1(G)$ and general Randi\'{c} index are depicted.

\begin{thm}\label{thm2}
Let $G$ be a graph with $n$ vertices and $m$ edges. Then for $\alpha \neq 1$, we have
\begin{eqnarray}
&(i)&I_{\mathscr{L}(\mathcal{Q})}^1(G)=1-\frac{1}{n^2}(n+2R_{-1}(G)),
\label{eq21}\\
&(ii)&I_\alpha^2(G)=\frac{1}{1-\alpha}\log\frac{M_\alpha^*}{n^\alpha}
\label{eq22}\\
&(iii) &I_\alpha^3(G)=\frac{1}{2^{1-\alpha}-1}
\Bigg(\frac{M_\alpha^*}{n^\alpha}-1\Bigg).\label{eq23}
\end{eqnarray}
where $R_{-1}(G)$ denotes the general Randi\'{c} index $R_\beta(G)$ of $G$ with $\beta=-1$ and $M_\alpha^*=\sum\limits_{i=1}^n|q_i|^\alpha$.
\end{thm}

\pf By substituting $\sum\limits_{i=1}^n \mu_i =\sum\limits_{i=1}^n q_i=n$ and $\sum\limits_{i=1}^n \mu_i^2
=\sum\limits_{i=1}^n q_i^2=n+2 \sum\limits_{i\sim j}\frac{1}{d_id_j}$ into equality (\ref{eq1}), we have
\begin{align*}
I^1_\mathscr{L}(G)&=1-\frac{1}{\left(\sum\limits_{j=1}^n|\mu_j|\right)^2}
\sum\limits_{i=1}^{n}|\mu_i|^2\\
&=1-\frac{1}{n^2}(n+2\sum\limits_{i\sim j}\frac{1}{d_id_j})\\
&=1-\frac{1}{n^2}(n+2R_{-1}(G)).
\end{align*}
\begin{align*}
I^1_\mathcal{Q}(G)&=1-\frac{1}{\left(\sum\limits_{j=1}^n|q_j|\right)^2}
\sum\limits_{i=1}^{n}|q_i|^2\\
&=1-\frac{1}{n^2}(n+2\sum\limits_{i\sim j}\frac{1}{d_id_j})\\
&=1-\frac{1}{n^2}(n+2R_{-1}(G)).
\end{align*}
The other two equalities can be obtained by substituting $\sum\limits_{i=1}^n \mu_i =\sum\limits_{i=1}^n q_i =2m$ and $M_\alpha^*=\sum\limits_{i=1}^n|q_i|^\alpha$ into equalities (\ref{eq2}) and (\ref{eq3}), respectively.\qed

From equality (\ref{eq21}), we can easily infer the relation of $I_{\mathscr{L}(\mathcal {Q})}^1(G)$ and the general Randi\'{c} index $R_{-1}(G)$ of $G$. And it implies that for a graph $G$, each upper (lower) bound of the general Randi\'{c} index $R_{-1}(G)$ of $G$ can deduce a lower (an upper) bound of $I_{\mathscr{L}(\mathcal {Q})}^1(G)$ directly. Easily to check the following extremal properties of $I_{\mathscr{L}(\mathcal {Q})}^1(G)$ employing the bounds in \cite{S,LY}.

\begin{cor}
i. For a graph $G$ with $n$ vertices and $m$ edges, if $n$ is odd, then we have
$$1-\frac{2}{n}+\frac{1}{n^2}\leq I_{\mathscr{L}(\mathcal {Q})}^1(G) \leq 1-\frac{1}{n-1},$$
if $n$ is even, then we have
$$1-\frac{2}{n}\leq I_{\mathscr{L}(\mathcal {Q})}^1(G) \leq 1-\frac{1}{n-1}.$$
with right equality if and only if $G$ is a complete graph, and with left equality if and only if $G$ is the disjoint union of $\frac{n}{2}$ paths of length 1 for $n$ is even, and is the disjoint union of $\frac{n-3}{2}$ paths of length 1 and a path of length 2 for $n$ is odd.

ii. Let $G$ be a graph with $n$ vertices and $m$ edges. The minimum degree of $G$ is $\delta$ and
the maximum degree of $G$ is $\Delta$. Then
$$1-\frac{1}{n}-\frac{1}{n\delta}\leq I_{\mathscr{L}(\mathcal {Q})}^1(G)\leq 1-\frac{1}{n}-\frac{1}{n\Delta}.$$
Equality occurs in both bounds if and only if $G$ is a regular graph.
\end{cor}\qed

3. Let $I(G)$ be the incidence matrix of graph $G$. For a graph $G$ with vertex set $V(G)=\{v_1, v_2, \ldots, v_n\}$ and edge set $E(G)=\{e_1, e_2, \ldots, e_m\}$, the $(i,j)$-entry of $I(G)$ is 1 if the vertex $v_i$ is incident with the edge $e_j$, and is 0 otherwise. As we know, $Q(G)=D(G)+A(G)=I(G)\cdot I^T(G)$. And if the eigenvalues of $Q(G)$ are $q_1, q_2, \ldots, q_n$, then $\sqrt{q_1}, \sqrt{q_2}, \ldots, \sqrt{q_n}$ are the singular values of $I(G)$. In addition, the incidence energy of graph $G$ is defined as $IE(G)=\sum\limits_{i=1}^n \sqrt{q_i}$. Similarly, we consider the connection between the generalized graph entropy $I_I^1(G)$ and the incidence energy $IE(G)$. And we get the following result.

\begin{thm}\label{thm3}
Let $G$ be a graph with $n$ vertices and $m$ edges. Then for $\alpha \neq 1$, we have
\begin{eqnarray}
&(i)& I_{I}^1(G)=1-\frac{2m}{IE^2(G)},\label{eq31}\\
&(ii)&I_\alpha^2(G)=\frac{1}{1-\alpha}\log\frac{M_\alpha^*}{IE^\alpha(G)},
\label{eq32}\\
&(iii) &I_\alpha^3(G)=\frac{1}{2^{1-\alpha}-1}
\Bigg(\frac{M_\alpha^*}{IE^\alpha(G)}-1\Bigg).\label{eq33}
\end{eqnarray}
where $IE(G)$ denotes the incidence energy of $G$ and $M_\alpha^*=\sum\limits_{i=1}^n(\sqrt{q_i})^\alpha$.
\end{thm}

\pf By substituting $IE(G)=\sum\limits_{i=1}^n \sqrt{q_i}$ and $\sum_{i=1}^n  q_i=tr(Q(G))=2m$ into equality (\ref{eq1}), we have
\begin{align*}
I^1_I(G)&=1-\frac{1}{\left(\sum\limits_{j=1}^n|\sqrt{q_j}|\right)^2}
\sum\limits_{i=1}^{n}|\sqrt{q_i}|^2\\
&=1-\frac{1}{IE^2(G)}\sum\limits_{i=1}^nq_i\\
&=1-\frac{2m}{IE^2(G)}.
\end{align*}

The other two equalities can be obtained by substituting $M_\alpha^*=\sum\limits_{i=1}^n\sqrt{q_i}^\alpha$ and $IE(G)=\sum\limits_{i=1}^n \sqrt{q_i}$ into equalities (\ref{eq2}) and (\ref{eq3}), respectively.\qed

Equality (\ref{eq31}) suggests the fact that for a graph $G$, each upper (lower) bound of the incidence energy $IE(G)$ of $G$ can deduce an upper (a lower) bound of $I_I^1(G)$. And applying some known bounds in \cite{GKM,JKM}, we can obtain the following extremal properties of the general graph entropy $I_{I}^1(G)$.

\begin{cor}
i. For a graph $G$ with $n$ vertices and $m$ edges, we have
$$ 0 \leq I_I^1(G) \leq 1-\frac{1}{n}.$$
The left equality holds if
and only if $m\leq 1$, whereas the right equality holds if and only if $m = 0$.

ii. Let $T$ be a tree of order $n$. Then we have
$$I_I^1(S_n) \leq I_I^1(T) \leq I_I^1(P_n),$$
where $S_n$ and $P_n$ denote the star and path of order $n$, respectively.
\end{cor}\qed

4. Let the graph $G$ be a connected graph whose vertices are $v_1, v_2, \ldots, v_n $. The distance matrix of $G$ is defined as $D(G) = [d_{ij}]$, where $d_{ij}$ is the distance between the vertices $v_i$ and $v_j$ in $G$. We denote the eigenvalues of $D(G)$ by $\mu_1, \mu_2, \ldots, \mu_n$. It is easy to verify that $\sum\limits_{i=1}^n\mu_i=0$ and $\sum\limits_{i=1}^n \mu_i^2=2\sum\limits_{1 \leq i <j \leq n}(d_{ij})^2$. The distance energy of the graph $G$ is $DE(G)=\sum\limits_{i=1}^n|\mu_i|$. And the $k$-th distance moment of $G$ is defined as $W_k(G)=\frac{1}{2}\sum\limits_{1 \leq i <j \leq n}(d_{ij})^k$. Particularly, $W(G)=W_1(G)$ and $WW(G)=\frac{1}{2}(W_2(G)+W_1(G))$, where $W(G)$ and $WW(G)$ respectively denote the Wiener index and hyper-Wiener index of $G$. We get the equality $W_2(G)=\frac{1}{2}\sum\limits_{1 \leq i <j \leq n}(d_{ij})^2=2WW(G)-W(G)$ by simple calculations. And the following theorem describes the equality relationship of the generalized graph entropy $I_D^1(G)$, $DE(G)$, $W(G)$, $WW(G)$ and so on.

\begin{thm}\label{thm4}
Let $G$ be a graph with $n$ vertices and $m$ edges. Then for $\alpha \neq 1$, we have
\begin{eqnarray}
&(i) &I_{D}^1(G)=1-\frac{4}{DE^2(G)}(2WW(G)-W(G)),\label{eq41}\\
&(ii)&I_\alpha^2(G)=\frac{1}{1-\alpha}\log\frac{M_\alpha^*}{DE^\alpha(G)},
\label{eq42}\\
&(iii)&I_\alpha^3(G)=\frac{1}{2^{1-\alpha}-1}
\Bigg(\frac{M_\alpha^*}{DE^\alpha(G)}-1\Bigg).\label{eq43}
\end{eqnarray}
where $M_\alpha^*=\sum\limits_{i=1}^n|\mu_i|^\alpha$ and $DE(G)$ denotes the distance energy of $G$. Here $W(G)$ and $WW(G)$ are the Wiener index and hyper-Wiener index of $G$.
\end{thm}

\pf By substituting $DE(G)=\sum\limits_{i=1}^n |\mu_i|$ and $W_2(G)=2WW(G)-W(G)$ into equality (\ref{eq1}), we have
\begin{align*}
I^1_D(G)&=1-\frac{1}{\left(\sum\limits_{j=1}^n|\mu_j|\right)^2}
\sum\limits_{i=1}^{n}|\mu_i|^2\\
&=1-\frac{1}{DE^2(G)}\sum\limits_{i=1}^n\mu_i^2\\
&=1-\frac{2}{DE^2(G)}\sum\limits_{1 \leq i <j \leq n}(d_{ij})^2\\
&=1-\frac{4W_2(G)}{DE^2(G)}\\
&=1-\frac{4}{DE^2(G)}(2WW(G)-W(G)).
\end{align*}

The other two equalities can be obtained by substituting $M_\alpha^*=\sum\limits_{i=1}^n|\mu_i|^\alpha$ and $DE(G)=\sum\limits_{i=1}^n |\mu_i|$ into equalities (\ref{eq2}) and (\ref{eq3}), respectively.\qed

From equality (\ref{eq41}), we can easily infer the relation of $I_D^1(G)$ and the distance energy $DE(G)$ of $G$, the Wiener index and the hyper-Wiener index of $G$. Together with some known bounds in \cite{RRGRAW}, we have the following corollary on $I_D^1(G)$.

\begin{cor}For a graph with $n$ vertices and $m$ edges, we have
$$ 0 \leq I_D^1(G) \leq 1-\frac{1}{n} .$$
\end{cor}\qed

5. Let $G$ be a simple undirected graph, and $G^\sigma$ be an oriented graph of $G$ with the orientation $\sigma$. The skew adjacency matrix of $G^\sigma$ is the $n \times n$ matrix $S(G^\sigma)=[s_{ij}]$, where $s_{ij}=1$ and $s_{ji}=-1$ if $\langle v_i,v_j\rangle$ is an arc of $G^\sigma$, otherwise $s_{ij}=s_{ji}=0$. Let $\lambda_1, \lambda_2, \ldots, \lambda_n$ be the eigenvalues of it. Obviously, $\lambda_1, \lambda_2, \ldots, \lambda_n$ are all pure images and also $\sum\limits_{i=1}^n \lambda_i=0$, $\sum\limits_{i=1}^n \lambda_i^2=-2m$. Then $\sum\limits_{i=1}^n|\lambda_i|^2=2m$. The skew energy of $G^\sigma$ is $SE(G^\sigma)=\sum\limits_{i=1}^n|\lambda_i|$. Now, we focus on the extremal properties of the generalized graph entropies $I_S^1(G)$, $I_\alpha^2(G)$ and $I_\alpha^3(G)$.

\begin{thm}\label{thm5}
Let $G^\sigma$ be an oriented graph with $n$ vertices and $m$ arcs. Then for $\alpha \neq 1$, we have
\begin{eqnarray}
&(i) &I_{S}^1(G^\sigma)=1-\frac{2m}{SE^2(G^\sigma)},\label{eq51}\\
&(ii)&I_\alpha^2(G^\sigma)=\frac{1}{1-\alpha}
\log\frac{M_\alpha^*}{SE^\alpha(G^\sigma)},\label{eq52}\\
&(iii)&I_\alpha^3(G^\sigma)=\frac{1}{2^{1-\alpha}-1}
\Bigg(\frac{M_\alpha^*}{SE^\alpha(G^\sigma)}-1\Bigg).\label{eq53}
\end{eqnarray}
where $SE(G^\sigma)$ denotes the skew energy of $G^\sigma$ and $M_\alpha^*=\sum\limits_{i=1}^n|\lambda_i|^\alpha$.
\end{thm}

\pf By substituting $SE(G^\sigma)=\sum\limits_{i=1}^n |\lambda_i|$ and $\sum\limits_{i=1}^n  |\lambda_i|^2=2m$ into equality (\ref{eq1}), we have
\begin{align*}
I^1_S(G^\sigma)&=1-\frac{1}{\left(\sum\limits_{j=1}^n|\lambda_j|\right)^2}
\sum\limits_{i=1}^{n}|\lambda_i|^2\\
&=1-\frac{2m}{SE^2(G^\sigma)}.
\end{align*}

The other two equalities can be obtained by substituting $M_\alpha^*=\sum\limits_{i=1}^n|\lambda_i|^\alpha$ and $SE(G^\sigma)=\sum\limits_{i=1}^n |\lambda_i|$ into equalities (\ref{eq2}) and (\ref{eq3}), respectively.\qed

Obviously, the equality (\ref{eq51}) implies the fact that for an oriented graph $G^\sigma$, each upper (lower) bound of the skew energy $SE(G^\sigma)$ of $G^\sigma$ can deduce an upper (a lower) bound of $I_S^1(G^\sigma)$. It is easy to check the following results based on some bounds in \cite{ABS}.

\begin{cor}
i. For an oriented graph $G^\sigma$ with $n$ vertices, $m$ arcs and maximum degree $\Delta$, we have
$$1-\frac{2m}{2m+n(n-1)|det(S(G^\sigma))|^{\frac{2}{n}}} \leq I_S^1(G^\sigma)\leq 1-\frac{1}{n} \leq 1-\frac{2m}{n^2\Delta} .$$
ii. Let $T^\sigma$ be an oriented tree of order $n$. We have
$$I_S^1(S_n^\sigma) \leq I_S^1(T^\sigma) \leq I_S^1(P_n^\sigma),$$
where $S_n^\sigma$ and $P_n^\sigma$ denote an oriented star and an oriented path with any orientation of order $n$, respectively. Equality holds if and only if the underlying tree $T_n$ satisfies that $T_n\cong S_n$ or $T_n \cong P_n$.
\end{cor}\qed

6. Let $G$ be a simple graph. The Randi\'{c} adjacency matrix of $G$ is defined as $R(G)=[r_{ij}]$, where $r_{ij}=\frac{1}{\sqrt{d_id_j}}$ if $v_i$ and $v_j$ are adjacent vertices of $G$, otherwise $r_{ij}=0$. Denote $\rho_1, \rho_2, \ldots, \rho_n$ be its eigenvalues. Obviously, $\sum\limits_{i=1}^n\rho_i=0$ and $\sum\limits_{i=1}^n\rho_i^2=tr(R^2(G))=2\sum\limits_{i\sim j}\frac{1}{d_id_j}$. The Randi\'{c} energy of the graph $G$ is defined as $RE(G)=\sum\limits_{i=1}^n|\rho_i|$.

\begin{thm}\label{thm6}
Let $G$ be a graph with $n$ vertices and $m$ edges. Then for $\alpha \neq 1$, we have
\begin{eqnarray}
&(i) &I_{R}^1(G)=1-\frac{2}{RE^2(G)}R_{-1}(G),\label{eq61}\\
&(ii)& I_\alpha^2(G)=\frac{1}{1-\alpha}\log\frac{M_\alpha^*}{RE^\alpha(G)},
\label{eq62}\\
&(iii)&I_\alpha^3(G)=\frac{1}{2^{1-\alpha}-1}
\Bigg(\frac{M_\alpha^*}{RE^\alpha(G)}-1\Bigg).\label{eq63}
\end{eqnarray}
where $RE(G)$ denotes the Randi\'{c} energy of $G$, and $R_{-1}(G)$ denotes the general Randi\'{c} index $R_{\beta}(G)$ of $G$ with $\beta=-1$ and $M_\alpha^*=\sum\limits_{i=1}^n|\rho_i|^\alpha$.
\end{thm}

\pf By substituting $RE(G)=\sum\limits_{i=1}^n |\rho_i|$ and $R_{-1}(G)=\sum\limits_{i\sim j}\frac{1}{d_id_j}$ into equality (\ref{eq1}), we have

\begin{align*}
I^1_R(G)&=1-\frac{1}{\left(\sum\limits_{j=1}^n|\rho_j|\right)^2}
\sum\limits_{i=1}^{n}|\rho_i|^2\\
&=1-\frac{1}{RE^2(G)}\sum\limits_{i=1}^n\rho_i^2\\
&=1-\frac{2}{RE^2(G)}\sum\limits_{i\sim j}\frac{1}{d_id_j}\\
&=1-\frac{2}{RE^2(G)}R_{-1}(G).
\end{align*}

The other two equalities can be obtained by substituting $M_\alpha^*=\sum\limits_{i=1}^n|\rho_i|^\alpha$ and $RE(G)=\sum\limits_{i=1}^n |\rho_i|$ into equalities (\ref{eq2}) and (\ref{eq3}), respectively.\qed

Equality (\ref{eq61}) suggests the relation of $I_R^1(G)$, the Randi\'{c} energy $RE(G)$ of $G$ and the general Randi\'{c} index $R_{-1}(G)$ of $G$. Together with some known bounds in \cite{BGGC}, we have the following corollary on $I_R^1(G)$.

\begin{cor}
For a graph with $n$ vertices and $m$ edges, we have
$$ I_R^1(G) \leq 1-\frac{1}{n} .$$
Equality is attained if and only if $G$ is the graph without edges, or if all its
vertices have degree one.
\end{cor}\qed

7. Let $G$ be a simple graph with vertex set $V(G)=\{v_1, v_2, \ldots, v_n\}$ and edge set $E(G)= \{e_1, e_2, \ldots, e_m\}$, and let $d_i$ be the degree of vertex
$v_i, i = 1, 2, \ldots, n$. Define an $n \times m$ matrix whose $(i, j)$-entry is $(d_i)^{-\frac{1}{2}}$ if $v_i$ is incident to $e_j$ and 0 otherwise. We call it the Randi\'{c} incidence matrix of $G$ and denote it by $I_R(G)$. Obviously, $I_R(G)=D(G)^{-\frac{1}{2}}I(G)$. Let $\sigma_1, \sigma_2, \ldots, \sigma_n$ be its singular values. And also $\sum\limits_{i=1}^n\sigma_i$ are defined as the Randi\'{c} incidence energy $I_RE(G)$ of the graph $G$. Let $U$ be the set of isolated vertices of $G$ and $W=V(G)-U$. Set $r=|W|$. Then we have $\sum\limits_{i=1}^n\sigma_i^2=r$. Particularly, $\sum\limits_{i=1}^n \sigma_i^2=n$ if $G$ has no isolated vertices.

\begin{thm}\label{thm7}
Let $G$ be a graph with $n$ vertices and $m$ edges. Let $U$ be the set of isolated vertices of $G$ and $W=V(G)-U$. Set $r=|W|$. Then for $\alpha \neq 1$, we have
\begin{eqnarray}
&(i) & I_{I_R}^1(G)=1-\frac{r}{I_RE^2(G)},\label{eq71}\\
&(ii)& I_\alpha^2(G)=\frac{1}{1-\alpha}\log\frac{M_\alpha^*}{I_RE^\alpha(G)},
\label{eq72}\\
&(iii)& I_\alpha^3(G)=\frac{1}{2^{1-\alpha}-1}
\Bigg(\frac{M_\alpha^*}{I_RE^\alpha(G)}-1\Bigg).\label{eq73}
\end{eqnarray}
where $I_RE(G)$ denotes the Randi\'{c} incidence energy of $G$ and $M_\alpha^*=\sum\limits_{i=1}^n|\sigma_i|^\alpha$.
\end{thm}

\pf By substituting $I_RE(G)=\sum\limits_{i=1}^n |\sigma_i|$ and $\sum\limits_{i=1}^n  |\sigma_i|^2=r$ into equality (\ref{eq1}), we have
\begin{align*}
I^1_{I_R}(G)&=1-\frac{1}{\left(\sum\limits_{j=1}^n|\sigma_j|\right)^2}
\sum\limits_{i=1}^{n}|\sigma_i|^2\\
&=1-\frac{r}{I_RE^2(G)}.
\end{align*}

The other two equalities can be obtained by substituting $M_\alpha^*=\sum\limits_{i=1}^n|\sigma_i|^\alpha$ and $I_RE(G)=\sum\limits_{i=1}^n |\sigma_i|$ into equalities (\ref{eq2}) and (\ref{eq3}), respectively.\qed

From equality (\ref{eq71}), we can easily infer the relation of $I_{I_R}^1(G)$ and the Randi\'{c} incidence energy $I_RE(G)$ of $G$. This equality tells us that for a graph $G$, each upper (lower) bound of the skew energy $I_RE(G)$ of $G$ can deduce an upper (a lower) bound of $I_{I_R}^1(G)$. Applying some known bounds in \cite{GHL1}, we can obtain the following extremal properties of the generalized graph entropy $I_{I_R}^1(G)$.

\begin{cor}
i. For a graph $G$ with $n$ vertices and $m$ edges, we have
$$I_{I_R}^1(G)\geq 1-\frac{r}{n},$$
the equality holds if and only if $G \cong K_2$.

ii. Let $G$ be a graph with $n$ vertices and $m$ edges. Then
$$I_{I_R}^1(G)\leq 1-\frac{r}{n^2-3n+4+2\sqrt{2(n-1)(n-2)}},$$
the equality holds if and only if $G \cong K_n$.

iii. Let $T$ be a tree of order $n$. We have
$$I_{I_R}^1(T) \leq I_{I_R}^1(S_n),$$
where $S_n$ denotes the star graph of order $n$.
\end{cor}\qed

8. Let $R_\beta(G)$ be the general Randi\'{c} matrix of graph $G$. Define $R_\beta(G)=[r_{ij}]$, where $r_{ij}=\frac{1}{(d_id_j)^\beta}$ if $v_i$ and $v_j$ are adjacent vertices of $G$, otherwise $r_{ij}=0$. Set $\gamma_1, \gamma_2, \ldots, \gamma_n$ be the eigenvalues of $R_\beta(G)$. By the definition of $R_\beta(G)$ we can get $R_\beta(G)=D(G)^\beta A(G)D(G)^\beta$ and $\sum\limits_{i=1}^n\gamma_i^2=tr(R_\beta^2(G))=2\sum\limits_{i\sim j}(d_id_j)^{2\beta}$ directly. The general Randi\'{c} energy is defined as $RE_\beta(G)=\sum\limits_{i=1}^n|\gamma_i|$. Similarly, we obtain the theorem as follows.

\begin{thm}\label{thm8}
Let $G$ be a graph with $n$ vertices and $m$ edges. Then for $\alpha \neq 1$, we have
\begin{eqnarray}
&(i)&I_{R_\beta}^1(G)=1-\frac{2}{RE_\beta^2(G)}R_{2\beta}(G),
\label{eq81}\\
&(ii)&I_\alpha^2(G)=\frac{1}{1-\alpha}
\log\frac{M_\alpha^*}{RE_\beta^\alpha(G)},\label{eq82}\\
&(iii)& I_\alpha^3(G)=\frac{1}{2^{1-\alpha}-1}
\Bigg(\frac{M_\alpha^*}{RE_\beta^\alpha(G)}-1\Bigg).\label{eq83}
\end{eqnarray}
where $RE_\beta(G)$ denotes the general Randi\'{c} energy of $G$, and $R_{2\beta}(G)$ denotes the general Randi\'{c} index of $G$ and $M_\alpha^*=\sum\limits_{i=1}^n|\gamma_i|^\alpha$.
\end{thm}

\pf By substituting $RE_\beta(G)=\sum\limits_{i=1}^n |\gamma_i|$ and $\sum\limits_{i=1}^n |\gamma_i|^2=2\sum\limits_{i\sim j}(d_id_j)^{2\beta}$ into equality (\ref{eq1}), we have
\begin{align*}
I^1_{R_\beta}(G)&=1-\frac{1}{\left(\sum\limits_{j=1}^n|\gamma_j|\right)^2}
\sum\limits_{i=1}^{n}|\gamma_i|^2\\
&=1-\frac{2}{RE_\beta^2(G)}\sum\limits_{i\sim j}(d_id_j)^{2\beta}\\
&=1-\frac{2}{RE_\beta^2(G)}R_{2\beta}(G).
\end{align*}

The other two equalities can be obtained by substituting $M_\alpha^*=\sum\limits_{i=1}^n|\gamma_i|^\alpha$ and $RE_\beta(G)=\sum\limits_{i=1}^n |\gamma_i|$ into equalities (\ref{eq2}) and (\ref{eq3}), respectively.\qed

9. Let $G$ be a simple undirected graph, and $G^\sigma$ be an oriented graph of $G$ with the orientation $\sigma$. The skew Randi\'{c} matrix of $G^\sigma$ is the $n \times n$ matrix $R_s(G^\sigma)=[(r_s)_{ij}]$, where $(r_s)_{ij}=(d_id_j)^{-\frac{1}{2}}$ and $(r_s)_{ji}=-(d_id_j)^{-\frac{1}{2}}$ if $\langle v_i,v_j\rangle$ is an arc of $G^\sigma$, otherwise $(r_s)_{ij}=(r_s)_{ji}=0$. Let $\rho_1, \rho_2, \ldots, \rho_n$ be the eigenvalues of $R_s(G^\sigma)$. It follows that $R_s(G^\sigma)=D(G)^{-\frac{1}{2}}S(G^\sigma)D(G)^{-\frac{1}{2}}$ and $\sum\limits_{i=1}^n\rho_i^2=tr(R_s^2(G^\sigma))=-2\sum\limits_{i\sim j}\frac{1}{d_id_j}=-2R_{-1}(G)$, which implies that $\sum\limits_{i=1}^n|\rho_i|^2=2R_{-1}(G)$. The skew Randi\'{c} energy is $RE_s(G^\sigma)=\sum\limits_{i=1}^n|\rho_i|$. We will give the accurate relation among $RE_s(G^\sigma)$, $I_{R_s}^1(G^\sigma)$, $I_{\alpha}^2(G^\sigma)$ and $I_{\alpha}^3(G^\sigma)$.

\begin{thm}\label{thm9}
Let $G^\sigma$ be an oriented graph with $n$ vertices and $m$ arcs. Then for $\alpha \neq 1$, we have
\begin{eqnarray}
&(i)&I_{R_S}^1(G^\sigma)=1-\frac{2}{RE_S^2(G^\sigma)}R_{-1}(G),
\label{eq91}\\
&(ii)&I_\alpha^2(G^\sigma)=\frac{1}{1-\alpha}
\log\frac{M_\alpha^*}{RE_S^\alpha(G^\sigma)},\label{eq92}\\
&(iii)&I_\alpha^3(G^\sigma)=\frac{1}{2^{1-\alpha}-1}
\Bigg(\frac{M_\alpha^*}{RE_S^\alpha(G^\sigma)}-1\Bigg).\label{eq93}
\end{eqnarray}
where $RE_S(G^\sigma)$ denotes the skew Randi\'{c} energy of $G^\sigma$, and $R_{-1}(G)$ denotes the general Randi\'{c} index of the underlying graph $G$ with $\beta=-1$ and $M_\alpha^*=\sum\limits_{i=1}^n|\rho_i|^\alpha$.
\end{thm}

\pf By substituting $RE_S(G^\sigma)=\sum\limits_{i=1}^n |\rho_i|$ and $\sum\limits_{i=1}^n|\rho_i|^2=2R_{-1}(G)$ into equality (\ref{eq1}), we have

\begin{align*}
I^1_{R_s}(G^\sigma)&=1-\frac{1}{\left(\sum\limits_{j=1}^n|\rho_j|\right)^2}
\sum\limits_{i=1}^{n}|\rho_i|^2\\
&=1-\frac{2}{RE_s^2(G)}R_{-1}(G).
\end{align*}

The other two equalities can be obtained by substituting $M_\alpha^*=\sum\limits_{i=1}^n|\rho_i|^\alpha$ and $RE_S(G^\sigma)=\sum\limits_{i=1}^n |\rho_i|$ into equalities (\ref{eq2}) and (\ref{eq3}), respectively.\qed

The equality (\ref{eq91}) states the relation of $I_{R_S}^1(G^\sigma)$ and the skew Randi\'{c} energy $RE_S(G^\sigma)$ of $G^\sigma$ and the general Randi\'{c} index $R_{-1}(G)$ of $G$. Together with some known bounds in \cite{GHL3}, we can obtain the following extremal properties of the generalized graph entropy $I_{R_S}^1(G)$.

\begin{cor}For an oriented graph $G^\sigma$ with $n$ vertices and $m$ arcs, we have
$$ I_{R_S}^1(G^\sigma) \leq 1-\frac{1}{n} .$$
\end{cor}\qed

For the above nine different entropies, we present the following results on implicit information inequality, which can be obtained by the method in \cite{DLS}.

\begin{thm}\label{thm10}
i. When $0<\alpha <1$, we have $I_\alpha^2 <I_\alpha^3\cdot \ln2$; and when $\alpha>1$, we have $I_\alpha^2> \frac{(1-2^{1-\alpha})\ln2}{\alpha-1}I_\alpha^3$.

ii. When $\alpha \geq 2$ and $0<\alpha <1$, we have $I_\alpha^3>I^1$; when $1<\alpha <2$, we have $I^1> (1-2^{1-\alpha})I_\alpha^3$.

iii. When $\alpha \geq 2$, we have
$$I_\alpha^2> \frac{(1-2^{1-\alpha})\ln2}{\alpha-1}I^1;$$
when $1<\alpha <2$, we have
$$I_\alpha^2> \frac{(1-2^{1-\alpha})^2\ln2}{\alpha-1}I^1;$$
when $0<\alpha <1$, we have $I_\alpha^2>I^1$.
\end{thm}

\section{Summary and Conclusion}\label{sec_conc}
In this paper, we proved interrelations between generalized graph entropies, distinct graph entropies and topological indices by means of inequalities. Graph energy and graph entropy are
well-defined concepts which have been introduced by Gutman \cite{Gutman_energy} and Mowshowitz \cite{M}, respectively. In terms of graph energy, various results have been obtained when proving
extremal results, see, e.g., \cite{ABS,BGGC,Gutman_energy,GHL1}. Also graph entropy is an important method introduced by Mowshowitz \cite{M} for determining the
structural information content of graphs that has been further developed by
many authors such as Bonchev \cite{B1}, K\"orner \cite{koerner_1973}, and Dehmer \cite{D}.

In view of the large amount of existing graph measures, the problem of deriving inequalities involving these measures has been only liitle investigated. For example, earlier work when proving
interrelations (inequalities) between graph entropies can be found in \cite{D,DK}. We argue that by proving such inequalities, we better understand the measures themselves and their behaviour. This may lead to novel
applications and to the dissemination of the results towards other disciplines.

\end{document}